\input amstex
\documentstyle{amsppt}
\magnification=\magstep1
\pagewidth{6.5truein}
\pageheight{9.0truein}

\def\rann{\operatorname{r{.}ann}}
\def\isosmd{\lesssim^{\oplus}}

\def\konealg{K_1^{\operatorname{alg}}}
\def\CC{{\Bbb C}}

\def\Ara{{\bf 1}}
\def\AGOPdiag{{\bf 2}}
\def\AGOPsep{{\bf 3}}
\def\Bla{{\bf 4}}
\def\Bro{{\bf 5}}
\def\BP{{\bf 6}}
\def\BPtwo{{\bf 7}}
\def\BPthree{{\bf 8}}
\def\CL{{\bf 9}}
\def\Coh{{\bf 10}}
\def\CJ{{\bf 11}}
\def\CR{{\bf 12}}
\def\Gnonsing{{\bf 13}}
\def\Gbull{{\bf 14}}
\def\Gvnrr{{\bf 15}}
\def\Gnonriesz{{\bf 16}}
\def\GW{{\bf 17}}
\def\HanI{{\bf 18}}
\def\HanII{{\bf 19}}
\def\Lin{{\bf 20}}
\def\LiRo{{\bf 21}}
\def\MM{{\bf 22}}
\def\Monc{{\bf 23}}
\def\Nic{{\bf 24}}
\def\Per{{\bf 25}}
\def\Rie{{\bf 26}}
\def\Ros{{\bf 27}}
\def\Sto{{\bf 28}}
\def\War{{\bf 29}}
\def\WO{{\bf 30}}
\def\Yu{{\bf 31}}
\def\Zone{{\bf 32}}
\def\Ztwo{{\bf 33}}
\def\Zthree{{\bf 34}}
\def\Zfour{{\bf 35}}

\topmatter

\pretitle{
\hbox{}
\vskip -0.9truein
\noindent {\sevenrm revised June 1999}
\vskip-0.0625truein

\noindent{\sevenrm to appear in Pacific
J.~Math.}
\vskip 0.9truein
         }
\title $K_1$ of separative exchange rings and C*-algebras
with real rank zero\endtitle

\rightheadtext{$K_1$ of separative exchange rings}

\author P. Ara, K. R. Goodearl, K. C. O'Meara, and R. Raphael
\endauthor

\address {P. Ara: Departament de Matem\`atiques, Universitat
Aut\`onoma de Bar\-ce\-lo\-na, 08193 Bellaterra (Barcelona),
Spain}\endaddress
\email para\@mat.uab.es\endemail

\address {K. R. Goodearl: Department of Mathematics,
University of California, Santa Barbara, California 93106,
USA}\endaddress
\email goodearl\@math.ucsb.edu\endemail

\address {K. C. O'Meara: Department of Mathematics, University
of  Canterbury, Christchurch, New Zealand}\endaddress
\email komeara\@math.canterbury.ac.nz\endemail

\address{R. Raphael: Department of Mathematics and
Statistics, Concordia University, Mont\-r\'eal, Qu\'ebec H4B 1R6,
Canada}\endaddress
\email raphael\@alcor.concordia.ca\endemail

\abstract For any (unital) exchange ring $R$ whose finitely
generated projective modules satisfy the separative
cancellation property ($A\oplus A\cong A\oplus B\cong
B\oplus B \implies A\cong B$), it is shown that all
invertible square matrices
over $R$ can be diagonalized by elementary row and column
operations. Consequently, the natural homomorphism $GL_1(R)
\rightarrow K_1(R)$ is surjective. In combination with a
result of Huaxin Lin, it follows that for any separative,
unital C*-algebra $A$ with real rank zero, the topological
$K_1(A)$ is naturally isomorphic to the unitary group
$U(A)$ modulo the connected component of the identity. This
verifies, in the separative case, a conjecture of Shuang
Zhang.
\endabstract

\thanks The research of the first author was partially
supported by grants from the DGICYT (Spain) and the Comissionat
per Universitats i Recerca de la Generalitat de Catalunya, that of
the second by a grant from the NSF (USA), and that of the fourth
by a grant from the NSERC (Canada).
\endthanks

\endtopmatter

\document

\head Introduction\endhead

The extent to which matrices over a ring $R$ can be
diagonalized is a measure of the complexity of $R$, as well
as a source of computational information about $R$ and its
free modules. Two natural properties offer themselves as
``best possible'': (1) that an arbitrary matrix can
be reduced to a diagonal matrix on left and right
multiplication by suitable invertible matrices, or (2) that
an arbitrary invertible matrix can be reduced to a diagonal
one by suitable elementary row and column operations. The
second property has an immediate K-theoretic benefit, in
that it implies that the Whitehead group $K_1(R)$ is a
natural quotient of the group of units of $R$. Our main
goal here is to prove property (2) for exchange rings
(definition below) satisfying a cancellation condition
which holds very widely (and conceivably for all exchange
rings). This theorem, when applied to C*-algebras with real
rank zero (also defined below), verifies a conjecture of
Shuang Zhang in an extensive class of C*-algebras.

The class of exchange rings has recently taken on a
unifying role for certain direct sum cancellation
problems in ring theory and operator algebra. In
particular, exchange rings encompass both (von Neumann)
regular rings (this is an old and easy observation) on the one
hand, and C*-algebras with real rank zero \cite{\AGOPsep,
Theorem 7.2} on the other. Within this class, a unifying
theme for a number of open problems is the property of
{\it separative cancellation\/} for finitely generated
projective modules, namely the condition
$$A\oplus A\cong A\oplus B\cong
B\oplus B \implies A\cong B$$
(see \cite{\AGOPdiag, \AGOPsep}). For example, if $R$ is a
separative exchange ring, then the (K-theoretic) stable
rank of $R$ can only be 1, 2, or $\infty$ \cite{\AGOPsep,
Theorem 3.3}, and every regular square matrix over $R$ is
equivalent (via multiplication by invertible matrices) to
a diagonal matrix \cite{\AGOPdiag, Theorem 2.4}. We prove
below that invertible matrices over separative exchange
rings can be diagonalized via elementary row and column
operations. Recently, Perera \cite{\Per} has applied our methods
to the problem of lifting units modulo an ideal $I$ in a ring $R$,
assuming that $I$ satisfies non-unital versions of separativity
and the exchange property. In this case, a unit $u$ of $R/I$ lifts
to a unit of $R$ if and only if the class of $u$ in $K_1(R/I)$ is
in the kernel of the connecting homomorphism $K_1(R/I) \rightarrow
K_0(I)$ \cite{\Per, Theorem 3.1}.

We defer discussion of the C*-algebraic aspects of our
results to Section 3, except for the following remark.
While earlier uses of the exchange property and
separativity for C*-algebras can easily be written out in
standard C*-theoretic terms -- e.g., with direct sums and
isomorphisms of finitely generated projective modules
replaced by orthogonal sums and Murray-von Neumann
equivalences of projections -- our present methods do not
lend themselves to such a translation. In particular,
although our main C*-algebraic application may be stated as
a diagonalization result for unitary matrices, all of the
steps in our proofs involve manipulations with non-unitary
matrices.

\medskip

Throughout the paper, we consider only unital rings and
C*-algebras. We reserve the term {\it elementary
operation\/} for the row (respectively, column) operation
in which a left (respectively, right) multiple of one row
(respectively, column) of a matrix is added to a different
row (respectively, column). Similarly, we reserve the name
{\it elementary matrix\/} for a transvection $I+re_{ij}$
where $I$ is an identity matrix, $e_{ij}$ is one of the
usual matrix units for some $i\ne j$, and $r$ is an element
of the base ring. Thus, as usual, an elementary row (respectively,
column) operation on a matrix $A$ corresponds to
multiplying $A$ on the left (respectively, right) by an
elementary matrix.

Note that while odd permutation matrices usually cannot be
expressed as products of elementary matrices, certain
signed permutation matrices can be. For example,
$$\left[ \matrix 0&1\\ -1&0 \endmatrix \right] = \left[
\matrix 1&1\\ 0&1 \endmatrix \right] \left[ \matrix 1&0\\
-1&1 \endmatrix \right] \left[ \matrix 1&1\\ 0&1 \endmatrix
\right] .$$
In particular, the operation of replacing rows $R_i$ and
$R_j$ (respectively, columns $C_i$ and $C_j$) with the rows
$R_j$ and $-R_i$ (respectively, the columns $C_j$ and
$-C_i$) can be achieved as a sequence of three elementary
operations. Therefore any entry of a matrix can be moved to
any other position by a sequence of elementary row and
column operations, at the possible expense of moving other entries
and multiplying some by $-1$.

For any ring $R$, let $E_n(R)$ denote the subgroup of
$GL_n(R)$ generated by the elementary matrices. If $GL_n(R)$
is generated by $E_n(R)$ together with the subgroup $D_n(R)$
of invertible
diagonal matrices, then $R$ is said to be a {\it
$GE_n$-ring\/}
\cite{\Coh, p\. 5}. Further, $R$ is a {\it $GE$-ring\/}
provided it is a
$GE_n$-ring for all $n$. If $R$ is a $GE_n$-ring
then $E_n(R)$ is a normal subgroup of $GL_n(R)$, and so
$GL_n(R)=D_n(R)E_n(R)= E_n(R)D_n(R)$. Of course, this means that
every invertible $n\times n$ matrix over a $GE_n$-ring can be
diagonalized using only elementary row (respectively, column)
operations.

It is easy to check
that all rings with stable rank $1$ are $GE$-rings. Note
that if
$R$ is a
$GE$-ring, then the natural homomorphism from $GL_1(R)$,
the group of units of $R$, to
$K_1(R)$ is surjective. For comparison, we recall the
well-known fact that if $R$ has stable rank $d$, then the
natural map
$GL_d(R) \rightarrow K_1(R)$ is surjective (e.g.,
\cite{\CR, Theorem 40.42}).

\head 1. Exchange rings and separativity \endhead

Although our notions and results will be right-left symmetric,
all modules considered in this paper will be right modules.
A module $M$ over a ring $R$ has the {\it finite exchange
property\/}
\cite{\CJ} if for every $R$-module $A$ and any decompositions 
$$A=M'\oplus N= A_1 \oplus \dots \oplus A_n$$  
with $M'\cong M$,
there exist submodules $A'_i\subseteq A_i$ such that
$$A=M'\oplus A'_1 \oplus \dots \oplus A'_n.$$  
(It
follows from the modular law that $A'_i$ must be a direct
summand of
$A_i$ for all $i$.) It should be emphasized that the direct
sums in this definition are
internal direct sums of submodules of $A$. One advantage of
the resulting internal direct sum decompositions (as
opposed to isomorphisms with external direct sums) rests on
the fact that direct summands with common complements are
isomorphic -- e.g., $N\cong
\bigoplus_{i=1}^n A'_i$ above since each of these summands
of $A$ has $M'$ as a complementary summand.

Following Warfield \cite{\War}, we say that $R$ is an {\it
exchange ring\/} if $R_R$ satisfies the finite exchange
property. By
\cite{\War, Corollary 2}, this definition is left-right
symmetric.  If $R$ is an exchange ring, then every finitely
generated projective
$R$-module has the finite exchange property (by \cite{\CJ,
Lemma 3.10}, the finite exchange property passes to finite
direct sums and to direct summands), and so the
endomorphism ring of any
 such module is an exchange ring. Further, idempotents lift
modulo all ideals of an exchange ring \cite{\Nic, Theorem
2.1, Corollary 1.3}.

The class of exchange rings is quite large. It includes all 
semiregular rings (i.e., rings which modulo the Jacobson
radical are (von Neumann) regular and have
idempotent-lifting),  all
$\pi$-regular rings, and more; see \cite{\Ara, \Sto, \War}.
Further, all unital C*-algebras with real rank zero are
exchange rings
\cite{\AGOPsep, Theorem 7.2}. 

The following criterion for exchange rings was obtained
independently by Nicholson and the second author.

\proclaim{Lemma 1.1} \cite{\GW, p\. 167; \Nic, Theorem 2.1} A
ring $R$ is an exchange ring if and only if for every element
$a\in R$ there exists an idempotent $e\in R$ such that $e\in
aR$ and $1-e\in (1-a)R$. \qed
\endproclaim

In the above lemma, it is equivalent to ask that for any
$a_1,a_2\in R$ with
$a_1R+a_2R=R$, there exists an idempotent
$e\in a_1R$ such that $1-e \in a_2R$. We shall also need the
analogous property corresponding to sums of more than two
right ideals:

\proclaim{Lemma 1.2} \cite{\Nic, Theorem 2.1, Proposition
1.11} Let $R$ be an exchange ring. If
$I_1,\dots,I_n$ are right ideals of $R$ such that
$I_1+\dots+I_n=R$, then there exist orthogonal idempotents
$e_1,\dots,e_n\in R$ such that
$e_1+\dots+e_n=1$ and $e_j \in I_j$ for all $j$.
\qed\endproclaim

We reiterate that a ring $R$ is {\it separative\/} provided the
following cancellation property holds for finitely generated
projective right (equivalently, left)
$R$-modules
$A$ and $B$:
$$A\oplus A\cong A\oplus B\cong B\oplus B \quad\implies\quad
A\cong B.$$ See \cite{\AGOPsep} for the origin of this
terminology and for a number of equivalent conditions. We
shall need the following one:

\proclaim{Lemma 1.3} \cite{\AGOPdiag, Proposition 1.2;
\AGOPsep, Lemma 2.1} A ring
$R$ is separative if and only if whenever $A,B,C$ are
finitely generated projective right $R$-modules such that
$A\oplus C\cong B\oplus C$ and
$C$ is isomorphic to direct summands of both $A^n$ and $B^n$
for some
$n$, then $A\cong B$. \qed\endproclaim

Note, in particular, that if $R$ is separative and $A,B,C$ are
finitely generated projective right $R$-modules, then we can
certainly cancel
$C$ from $A\oplus C\cong B\oplus C$ whenever $A$ and $B$
are generators in
$\text{Mod-}R$.

Separativity seems to hold quite widely within the class of
exchange rings; for instance, it holds for all known classes
of regular rings (cf\. \cite{\AGOPsep}). In fact, the existence of
non-separative exchange rings is an open problem.

It is clear from either form of the condition that a ring $R$
is separative in case the finitely generated projective
$R$-modules enjoy cancellation with respect to direct sums,
which in turn holds in case
$R$ has stable rank 1. In fact, for exchange rings,
cancellation of finitely generated projective modules is
equivalent to stable rank 1
\cite{\Yu, Theorem 9}. Separativity, however, is much weaker
than stable rank 1. For example, any regular right
self-injective ring is separative (e.g., \cite{\Gvnrr,
Theorem 10.34(b)}), but such rings can have infinite stable
rank -- e.g., the ring of all linear transformations on an
infinite dimensional vector space.

\head 2. $K_1$ of separative exchange rings\endhead

We use the notation $A\isosmd B$ to denote that a module $A$
is isomorphic to a direct summand of a module $B$.

\proclaim{Lemma 2.1} Let $R$ be an exchange ring and
$e_1,\dots,e_n\in R$ idempotents. Then there exists an
idempotent $e\in e_1R+\dots+e_nR$ such that $e_1R\le eR$ and
$e_iR\isosmd eR$ for all $i$. In particular,
$ReR= Re_1R +\dots+ Re_nR$.\endproclaim

\demo{Proof} By induction, it suffices to do the case $n=2$.
Now
$$R= e_1R\oplus (1-e_1)R= e_2R\oplus (1-e_2)R$$ and $e_1R$
has the finite exchange property, so there exist
decompositions $e_2R= A\oplus B$ and $(1-e_2)R= A'\oplus B'$
such that
$R= e_1R\oplus A\oplus A'$. Then we can choose an idempotent
$e\in R$ such that $eR= e_1R\oplus A$. Obviously $e_1R\le
eR$, and since
$$e_1R\cong R/(A\oplus A')\cong B\oplus B',$$ we have $e_2R
\isosmd A\oplus e_1R= eR$. \qed\enddemo

\proclaim{Corollary 2.2} Let $R$ be an exchange ring and
$a\in R$ such that $RaR=R$. Then there exist idempotents
$e\in aR$ and $f\in Ra$ such that $ReR= RfR= R$.\endproclaim

\demo{Proof} Write $R= \sum_{i=1}^n x_iaR$ for some $x_i$. By
Lemma 1.2, there exist orthogonal idempotents
$g_1,\dots,g_n\in R$ such that $g_1+\dots+g_n=1$ and $g_i\in
x_iaR$ for all $i$. Set $g_i=x_iay_i$ with $y_i= y_ig_i$.
Then $e_i := ay_ix_i$ is an idempotent in $aR$ and $e_iR\cong
g_iR$. By Lemma 2.1, there exists an idempotent $e\in
\sum_{i=1}^n e_iR$ such that $e_iR\isosmd eR$ for all
$i$. Then $e\in aR$ and $g_iR \isosmd eR$ for all $i$, so all
$g_i\in ReR$, and thus $ReR=R$.

The existence of $f$ follows by symmetry. \qed\enddemo

\proclaim{Lemma 2.3} Let $R$ be any ring and $A\in GL_n(R)$.
If $A$ has an idempotent entry, then $A$ can be reduced by
elementary row and column operations to the form
$$\left[ \matrix 1&0&\cdots&0\\ 0&*&\cdots&*\\
\vdots&\vdots&&\vdots\\ 0&*&\cdots&* \endmatrix \right]
.$$\endproclaim

\demo{Proof} By elementary operations, we can move the
idempotent entry, call it $e$, into the $1,1$ position. If
$n=1$, then $e$ is invertible, so $e=1$ and we are done. Now
assume that $n>1$, and let
$$\left[ \matrix e&b_2&b_3&\cdots&b_n \endmatrix \right]$$ be
the first row of $A$. By elementary column operations, we can
subtract $eb_i$ from the $i$-th entry for each $i\ge2$. Thus,
we can assume that $b_2,\dots,b_n \in (1-e)R$. Since $A$ is
invertible,
$eR+b_2R+ \dots+ b_nR=R$, and so it follows that $b_2R+
\dots+ b_nR= (1-e)R$. Hence, by elementary column operations
we can add $1-e$ to the first entry. Now we have a $1$ in the
$1,1$ position, and the rest is routine. \qed\enddemo

Since we shall need to perform a number of operations on the
top rows of invertible matrices, it is convenient to work
with the rows alone. Recall that any row $\left[
\matrix a_1&a_2&\cdots&a_n \endmatrix \right]$ of an
invertible matrix over a ring $R$ is {\it right
unimodular\/}, that is, $\sum_{i=1}^n a_iR =R$. Elementary
column operations apply to such a row just by viewing it as a
$1\times n$ matrix. Such operations amount to multiplying the
row on the right by an elementary matrix. Since our rings
need not be commutative, elementary column operations can
only introduce right-hand coefficients.

\proclaim{Lemma 2.4} Let $R$ be an exchange ring and $\alpha=
\left[
\matrix a_1&a_2&\cdots&a_n \endmatrix \right]$ a right
unimodular row over $R$. Then $\alpha$ can be transformed by
elementary column operations to a row $\left[ \matrix
b_1&b_2&\cdots&b_n
\endmatrix \right]$ such that $R= b_1R\oplus \cdots\oplus
b_nR$ and each $b_i\in a_iRa_i$. \endproclaim

\demo{Proof} Since $\sum_{i=1}^n a_iR =R$, Lemma 1.2 gives us
orthogonal idempotents $e_1,\dots,e_n\in R$ such that
$e_1+\dots+e_n=1$ and $e_i\in a_iR$ for all $i$, say $e_i=
a_ir_i$. By elementary column operations, we can subtract
$e_ia_1= a_ir_ia_1$ from the first entry of $\alpha$ for each
$i\ge2$. This transforms $\alpha$ to $\alpha'= \left[
\matrix e_1a_1&a_2&a_3&\cdots&a_n \endmatrix \right]$. Note
that $e_1\in e_1a_1R$. Thus, we can repeat the above process
for each entry, and transform $\alpha'$ to the row $\left[
\matrix e_1a_1&e_2a_2&\cdots&e_na_n \endmatrix \right]$, with
entries $e_ia_i\in a_iRa_i$. Moreover, $e_ia_iR= e_iR$, and
therefore
$R= \bigoplus_{i=1}^n e_ia_iR$. \qed\enddemo

\proclaim{Corollary 2.5} Let $R$ be an exchange ring and
$\alpha= \left[
\matrix a_1&a_2&\cdots&a_n \endmatrix \right]$ a right
unimodular row over $R$, with $n\ge2$. Then $\alpha$ can be
transformed by elementary column operations to a row $\left[
\matrix b_1&b_2&\cdots&b_n
\endmatrix \right]$ such that $Rb_1R= R$ and $b_i\in a_iRa_i$
for all $i\ge2$. \endproclaim

\demo{Proof} By Lemma 2.4, we may assume that $R=
\bigoplus_{i=1}^n a_iR$. It follows that all $a_i\in Rb_1$
where $b_1= a_1+\dots+a_n$ (multiply
$b_1$ on the left by the orthogonal idempotents arising from
the given decomposition of $R_R$). Thus $Rb_1R=R$. By
elementary column operations, we can add $a_2,\dots,a_n$ to
the first entry of $\alpha$, and thus transform it to $\left[
\matrix b_1&a_2&\cdots&a_n
\endmatrix \right]$. \qed\enddemo

Recall that an element $x$ in a ring $R$ is {\it (von Neumann)
regular\/} provided there exists an element $y\in R$ such
that $xyx=x$, equivalently, provided $xR$ is a direct summand
of $R_R$. If
$y$ can be chosen to be a unit in $R$, then $x$ is said to be
{\it unit-regular\/}. A regular element $x\in R$ is
unit-regular if and only if $R/xR \cong \rann(x)$, where
$\rann(x)$ denotes the right annihilator of $x$ in $R$ (cf\.
\cite{\Gvnrr, proof of Theorem 4.1}).

\proclaim{Corollary 2.6} Let $R$ be an exchange ring and
$\alpha= \left[
\matrix a_1&a_2&\cdots&a_n \endmatrix \right]$ a right
unimodular row over $R$, with $n\ge2$. Then $\alpha$ can be
transformed by elementary column operations to a row $\left[
\matrix c_1&c_2&\cdots&c_n
\endmatrix \right]$ such that $c_2$ is a regular element,
$c_2\in Ra_2$, and $c_2R= (1-g)R$ for an idempotent $g$ with
$RgR=R$.
\endproclaim

\demo{Proof} By Corollary 2.5, we may assume that $Ra_1R=R$.
By Corollary 2.2, there exists an idempotent $e\in a_1R$ such
that $ReR=R$. By elementary column operations, we can
subtract $ea_2$ from the second entry of $\alpha$, so there
is no loss of generality in assuming that
$a_2\in (1-e)R$. (At this stage, our current $a_2$ is only a
\underbar{left} multiple of the original $a_2$. This is why
the conclusions of the lemma state $c_2\in Ra_2$ rather than
$c_2 \in a_2Ra_2$.) Now using Lemma 2.4, we can transform
$\alpha$ to a row
$\left[
\matrix c_1&c_2&\cdots&c_n
\endmatrix \right]$ such that $R= \bigoplus_{i=1}^n c_iR$ and
$c_2\in a_2Ra_2$. Then $c_2R= (1-g)R$ for some idempotent
$g$, and $c_2$ is regular. Moreover, $(1-g)R= c_2R\subseteq
a_2R\subseteq (1-e)R$ and so $Re\subseteq Rg$. Therefore
$RgR=R$. \qed\enddemo

\proclaim{Lemma 2.7} Let $R$ be an exchange ring and $A\in
GL_n(R)$, with
$n\ge2$. Then $A$ can be transformed by elementary row and
column operations to a matrix whose $1,1$ entry $d$ is
regular, with $dR= (1-p)R$ and $Rd= R(1-q)$ for some
idempotents $p,q$ such that $RpR= RqR=R$. \endproclaim

\demo{Proof} By Corollary 2.6, we can assume that the $1,2$
entry of $A$ is a regular element $c$ such that $cR= (1-g)R$
for some idempotent $g$ with $RgR= R$. With elementary
operations, we can move $c$ to the
$2,1$ position.

Now apply the transpose of Corollary 2.6 to the first column
of $A$. Thus,
$A$ can be transformed by elementary row operations to a
matrix whose
$2,1$ entry is a regular element $d$ such that $d\in cR$ and
$Rd= R(1-q)$ for some idempotent $q$ with $RqR=R$. Since $d$
is regular, $dR= (1-p)R$ for some idempotent $p$. Then
$(1-p)R\subseteq (1-g)R$, whence
$Rg\subseteq Rp$ and so $RpR=R$.

Finally, use elementary operations to move $d$ to the $1,1$
position.
\qed\enddemo

\proclaim{Theorem 2.8} If $R$ is a separative exchange ring,
then $R$ is a GE-ring, and so the natural homomorphism
$GL_1(R) \rightarrow K_1(R)$ is surjective.
\endproclaim

\demo{Proof} We need to show that $R$ is a $GE_n$-ring for
all $n$. This is trivial for $n=1$, hence we assume, by induction,
that
$n\ge2$ and that $R$ is a $GE_{n-1}$-ring. Let $A$ be an
arbitrary invertible $n\times n$ matrix over $R$.

By Lemma 2.7, we may assume that the $1,1$ entry $d$ of $A$
is regular, with $dR= (1-p)R$ and $Rd= R(1-q)$ for some
idempotents $p,q$ such that $RpR= RqR=R$. We claim that $d$
is unit-regular. Note that because $RpR= RqR=R$, the
projective modules $pR$ and $qR$ are generators.

Now $R= \rann(d)\oplus B= dR\oplus C$ for some $B,C$, and we
have to prove that $\rann(d)\cong C$. Since $B\cong dB= dR$,
we have
$\rann(d)\oplus B\cong C\oplus B$. From $Rd=R(1-q)$, we get
$\rann(d)=qR$ and so $\rann(d)$ is a generator. Since $C\cong
R/dR\cong pR$, we see that $C$ is a generator too. By Lemma
1.3, $\rann(d)\cong C$ as desired.

The unit-regularity of $d$ gives $d=ue$ for some unit
$u$ and idempotent $e$. Set
$$U= \left[ \matrix u&0&\cdots&0\\ 0&1&\cdots&0\\
&&\ddots\\ 0&0&\cdots&1 \endmatrix \right] ;$$
then the matrix $U^{-1}A$ has an
idempotent entry. By Lemma 2.3, there exist $E,F\in
E_n(R)$ such that
$$EU^{-1}AF= \left[ \matrix 1&0\\ 0&A' \endmatrix \right]$$
where
$A'\in GL_{n-1}(R)$. By our induction hypothesis, $A'\in
E_{n-1}(R)D_{n-1}(R)$. It follows that
$$A\in D_n(R)E_n(R)D_n(R)E_n(R),$$
and therefore we have shown that $R$ is a $GE_n$-ring. This
establishes the induction step and completes the proof.
\qed\enddemo

\definition{Remarks 2.9} {\bf (a)} Observe that the proof of
Theorem 2.8 did not use the full force of separativity, only
the cancellation property
$(A\oplus C\cong B\oplus C \implies A\cong B)$ for
finitely generated projective $R$-modules $A,B,C$ with $A$
and $B$ generators. 

{\bf (b)} Theorem 2.8 includes, in particular, the result of
Menal and Moncasi that every factor ring of a right
self-injective ring is a $GE$-ring
\cite{\MM, Theorem 2.2}. To make the connection explicit,
recall that right self-injective rings are semiregular (e.g.,
\cite{\Gnonsing, Theorem 2.16, Lemma 2.18}) and hence
exchange; thus, all their factor rings are exchange rings.
Further, any right self-injective ring is separative (e.g.,
\cite{\Gbull, Theorem 3}). It follows that factor rings of
right self-injective rings
are separative
\cite{\AGOPsep, Theorem 4.2}.

{\bf (c)} As a special case of Theorem 2.8, we obtain that
any separative regular ring is a $GE$-ring, which gives a 
partial affirmative answer to a question of Moncasi
\cite{\Monc, Questi\'o 5}.

{\bf (d)} In the situation of Theorem 2.8, one naturally asks for
a description of the kernel of the epimorphism $GL_1(R)
\rightarrow K_1(R)$. This has been answered for unit-regular rings
and regular right-self-injective rings by Menal and Moncasi
\cite{\MM, Theorems 1.6, 2.6}, and for exchange rings with
primitive factors artinian by Chen and Li \cite{\CL, Theorem 3}.
In all the above cases, $K_1(R)\cong GL_1(R)^{\text{ab}}$ provided
$\frac12 \in R$ \cite{\MM, Theorems 1.7, 2.6; \CL, Corollary 7}.
Further, $K_1(R)\cong GL_1(R)^{\text{ab}}$ when $R$ is either a
C*-algebra with unitary 1-stable range or an AW*-algebra
\cite{\MM, Theorem 1.3, Corollary 2.11} (here the algebraic $K_1$
is meant). The unit-regular and AW* results correct and extend
earlier work of Handelman
\cite{\HanI, Theorem 2.4; \HanII, Theorem 7}.
\enddefinition

\proclaim{Theorem 2.10} If $R$ is a separative exchange ring
and $A$ is a (von Neumann) regular $n\times n$ matrix over
$R$, then $A$ can be diagonalized using elementary row and
column operations. \endproclaim

\demo{Proof} By \cite{\AGOPdiag, Theorem 2.4}, there exist
$P,Q \in GL_n(R)$ such that $PAQ$ is diagonal. By Theorem
2.8, $P=U_1V_1$ and $Q=V_2U_2$, where $U_1,U_2\in D_n(R)$
and $V_1,V_2\in E_n(R)$. So $V_1AV_2$ is a diagonal matrix obtained from $A$ 
by elementary row and column operations.
\qed\enddemo

\definition{Remark 2.11} When applying Theorem 2.10, note the
distinction between invertible matrices and general matrices. An
invertible matrix over a separative exchange ring can be
diagonalized from either side (by Theorem 2.8), whereas the
diagonalization of a general regular matrix sometimes requires
elementary operations on both the rows and the columns. For
example, the $2\times 2$ matrix $\left[ \smallmatrix 1&1\\0&0
\endsmallmatrix \right]$ over a field cannot be diagonalized
using only elementary row operations. 
\enddefinition

\definition{Example 2.12} Non-regular matrices over
separative exchange rings need not be diagonalizable by
elementary operations, even over finite dimensional algebras.
For example, choose a field $F$ and let
$$R= F[x_1,x_2,x_3,x_4] / \langle x_1,x_2,x_3,x_4 \rangle^2.$$
Then $R$ has a basis $1,a_1,a_2,a_3,a_4$ such that $a_ia_j=0$ for all $i,j$.
Since $R$ is clearly semiregular, it is an exchange ring;
separativity is an easy exercise. In fact, since $R$ is
artinian, it has stable rank $1$. Recall that this also
implies that $R$ is a $GE$-ring. Observe that every element
of
$R$ is a sum of a scalar plus a nilpotent element, and that
the product of any two nilpotent elements of $R$ is zero.

Now consider the matrix $A= \left[ \smallmatrix a_1&a_2\\
a_3&a_4
\endsmallmatrix \right]$, whose entries are linearly
independent nilpotent elements of $R$. We claim that any
sequence of elementary row or column operations on $A$ can
only produce a matrix whose entries are linearly independent
nilpotent elements. For instance, consider a product
$$\left[ \matrix 1&b\\ 0&1 \endmatrix \right] \left[ \matrix
c_{11}&c_{12}\\ c_{21}&c_{22} \endmatrix \right] = \left[
\matrix c_{11}+bc_{21}&c_{12}+bc_{22}\\ c_{21}&c_{22}
\endmatrix \right]$$ where $c_{11},c_{12},c_{21},c_{22}$ are
linearly independent nilpotent elements. Then $b=\beta+n$ for
some $\beta\in F$ and some nilpotent element $n$, whence
$bc_{21}= \beta c_{21}$ and $bc_{22}= \beta c_{22}$, and so
the entries in the matrix product above are linearly
independent (they are clearly nilpotent). The same thing
happens with other elementary operations, establishing the
claim.

Therefore no sequence of elementary operations on $A$ can
produce a matrix with a zero entry. In particular, $A$ cannot
be diagonalized by elementary operations. Since $R$ is a
$GE$-ring, it follows that $A$ cannot be diagonalized by
invertible matrices either, i.e., there do not exist
$X,Y\in GL_2(R)$ such that $XAY$ is diagonal. Thus the first of
the natural properties discussed in the introduction is not
implied by the second.
\qed\enddefinition

\head 3. $K_1$ of separative C*-algebras with real rank
zero\endhead

In connection with his work on the structure of
multiplier algebras (e.g., \cite{\Zone,
\Ztwo, \Zthree}), Shuang Zhang has conjectured [unpublished]
that if
$A$ is any unital C*-algebra with real rank zero, the
topological
$K_1(A)$ is isomorphic to the unitary group $U(A)$ modulo the
connected component of the identity, $U(A)^\circ$. We confirm this
conjecture in case $A$ is separative, which at the same time
provides a unified approach to all known cases of the
conjecture. The main interest of Zhang's conjecture is in the case
when the stable rank of $A$ is greater than $1$, since it has long
been known that $K^{\text{top}}_1(A) \cong U(A)/U(A)^\circ$
for all unital C*-algebras $A$ with stable rank $1$ (e.g., this is
equivalent to
\cite{\Rie, Theorem 2.10}).

We consider only unital, complex C*-algebras in this section,
and we refer the reader to \cite{\Bla, \WO} for 
 background and notation for C*-algebras. In particular, we
use $\sim$ and $\lesssim$ to denote Murray-von Neumann
equivalence and subequivalence of projections, and we write 
$M_\infty(A)$ for the (non-unital) algebra  of
$\omega\times\omega$ matrices
 with only finitely many nonzero entries from an
algebra $A$. We write
$U(A)$ for the unitary group of a unital C*-algebra $A$, and
$U(A)^\circ$ for the connected component of the identity in
$U(A)$. 

In the theory of operator algebras, it is customary to write
$K_1(A)$ for the topological $K_1$-group of $A$ (e.g.,
\cite{\Bla, Definition 8.1.1; \WO, Definition 7.1.1}), and we
shall follow that practice here. Thus, $K_1(A)= GL_\infty(A)/
GL_\infty(A)^\circ$. We then use the notation
$\konealg(A)$ to denote the algebraic $K_1$-group of $A$.
Since
$\konealg(A)$ is the abelianization of $GL_\infty(A)$ (e.g.,
\cite{\Ros, Proposition 2.1.4, Definition 2.1.5}) and $K_1(A)$
is abelian (e.g., \cite{\Bla, Proposition 8.1.3; \WO,
Proposition 7.1.2}), there is a natural surjective
homomorphism $\konealg(A)
\rightarrow K_1(A)$. Finally, following Brown
\cite{\Bro, p\. 116}, we say that
$A$ has {\it
$K_1$-surjectivity\/} (respectively, {\it
$K_1$-injectivity\/}) provided the natural homomorphism
$U(A)/U(A)^\circ \rightarrow K_1(A)$ is surjective
(respectively, injective).

The concept of {\it real rank zero\/} for a C*-algebra $A$
has a number of equivalent characterizations (see
\cite{\BP}). One is the requirement that each self-adjoint
element of $A$ can be approximated arbitrarily closely by
real linear combinations of orthogonal projections. (This is
usually phrased as saying that the set of self-adjoint
elements of $A$ with finite spectrum is dense in the set of
all self-adjoint elements.) It was proved in \cite{\AGOPsep,
Theorem 7.2} that $A$ has real rank zero if and only if it is
an exchange ring. Hence, the C*-algebras with real rank zero are
exactly the C*-algebras to which our results above can be applied.

Given a C*-algebra $A$, all idempotents in matrix algebras
$M_n(A)$ are equivalent to projections (e.g., \cite{\Bla,
Proposition 4.6.2; \Ros, Proposition 6.3.12}).  Hence, $A$ is
separative if and only if
$$p\oplus p\sim p\oplus q\sim q\oplus q \qquad\implies\qquad
p\sim q$$  for  projections $p,q\in M_\infty(A)$. An
equivalent condition (analogous to Lemma 1.3) is that
$p\oplus r\sim q\oplus r$ $\implies$ $p\sim q$ whenever
$r\lesssim n{.}p$ and $r\lesssim n{.}q$ for some $n$.
Separativity in $A$ is thus equivalent to the requirement that all
matrix algebras $M_n(A)$ satisfy the {\it weak cancellation\/}
introduced by Brown and Pedersen
\cite{\Bro, p\. 116; \BPtwo, p\. 114}. They have shown that every
extremally rich C*-algebra (see \cite{\BPtwo, p\. 125}) with
real rank zero is separative (\cite{\BPthree}, announced in
\cite{\Bro, p\. 116}). We would like to emphasize the question of
whether non-separative exchange rings exist by focusing on the C*
case:

\definition{Problem} Is every C*-algebra with real rank zero
separative?
\enddefinition

By combining Theorem 2.8 with a result of Lin, we obtain the
following theorem.

\proclaim{Theorem 3.1} If $A$ is a separative, unital
C*-algebra with real rank zero, then the natural map
$U(A)/U(A)^\circ \rightarrow K_1(A)$ is an isomorphism.
\endproclaim

\demo{Proof} Lin proved $K_1$-injectivity for C*-algebras
with real rank zero in
\cite{\Lin, Lemma 2.2}. Hence, it only remains to show
$K_1$-surjectivity. It is a standard fact that
$U(A)$ and $GL_1(A)$ have the same image in $K_1(A)$
(e.g.,
\cite{\Bla, pp\. 66, 67} or
\cite{\WO, proof of Proposition 4.2.6}). Now the natural map
$GL_1(A) \rightarrow K_1(A)$ factors as the composition of
natural maps
$GL_1(A) \rightarrow \konealg(A) \rightarrow K_1(A)$, the
second of which is surjective. Since $A$ has real rank zero,
it is an exchange ring, and so the map $GL_1(A) \rightarrow
\konealg(A)$ is surjective by Theorem 2.8. Therefore the
image of
$U(A)$ in $K_1(A)$ is all of $K_1(A)$, as desired.
\qed\enddemo

Brown and Pedersen have proved that every separative,
extremally rich C*-algebra has $K_1$-surjectivity
(\cite{\BPthree}, announced in
\cite{\Bro, p\. 116; \BPtwo, p\. 114}). Since there are
C*-algebras with real rank zero that are not extremally rich
\cite{\Bro, Example, p\. 117}, Theorem 3.1 can be viewed as a
partial extension of the Brown-Pedersen result within
the class of C*-algebras with real rank zero.

We thank the referee for the following remark.

\definition{Remark 3.2} While Theorem 3.1 is neither unexpected
nor new in the case of stable rank $1$ (cf\. the result of Rieffel
cited above), it is perhaps surprising that there are many
C*-algebras of real rank zero and stable rank $2$ to which the
theorem applies. To see this, consider C*-algebra extensions
$$0\rightarrow I\rightarrow A\rightarrow B\rightarrow 0$$
in which $I$ and $B$ have real rank zero and $A$ is unital. By
theorems of Zhang (\cite{\Zfour, Lemma 2.4}; cf\. \cite{\BP,
Theorem 3.14}) and Lin and R\o rdam  \cite{\LiRo, Proposition 4},
$A$ has real rank zero if and only if projections lift from $B$ to
$A$, if and only if the
connecting map $K_0(B) \rightarrow K_1(I)$ in topological K-theory
vanishes. In this case, by
\cite{\AGOPsep, Theorem 7.5}, $A$ will be separative provided $I$
and $B$ are both separative, and in particular if $I$ and $B$ have
stable rank $1$. However, by \cite{\LiRo, Proposition 4}, if $I$
and $B$ have stable rank $1$, then $A$ will have stable rank $2$
provided the connecting map
$K_1(B) \rightarrow K_0(I)$ does not vanish. It is easy to find
specific extensions satisfying the above conditions, such as the
examples analyzed in \cite{\LiRo, End of Section 1} or
\cite{\Gnonriesz}.
\enddefinition

We conclude with an application of Theorem 3.1 that extends
an argument of Brown
\cite{\Bro, Theorem 1}, relating homotopy and unitary
equivalence of projections, to a wider context within real
rank zero. Projections $p$ and
$q$ in a C*-algebra
$A$ are {\it unitarily equivalent\/} provided there exists a
unitary element
$u\in A$ such that $upu^* = q$; they are {\it homotopic\/}
provided there is a continuous path $f : [0,1] \rightarrow
\{\text{projections in\ } A\}$ such that $f(0)=p$
and
$f(1)=q$. It is a standard fact that homotopic projections are
unitarily equivalent (e.g.,
\cite{\Bla, Propositions 4.3.3, 4.6.5; \WO, Proposition
5.2.10}).

\proclaim{Theorem 3.3} Let $A$ be a separative, unital
C*-algebra with real rank zero, let $p,q\in A$ be
projections, and let $B=
\overline{ApA} +\CC\cdot 1$. Then $p$ and $q$ are homotopic
in $A$ if and only if $q\in B$ and $p,q$ are unitarily
equivalent in $B$.
\endproclaim

\demo{Proof} If $p$ and $q$ are homotopic in $A$, they are
connected by a path of projections within $A$. Each
projection along this path is homotopic to $p$ and hence is
unitarily equivalent to $p$. Thus, these projections all lie
in $ApA$. In particular, $q\in B$, and $p$ and $q$ are
homotopic in
$B$. Consequently, $p$ and $q$ must be unitarily equivalent
in $B$.

Conversely, assume that $q\in B$ and $p,q$ are unitarily
equivalent in
$B$. By \cite{\BP, Corollary 2.8, Theorem 2.5}, the closed
ideal $I=
\overline{ApA}$ has real rank zero (as a non-unital
C*-algebra), and so the unital C*-algebras
$B$ and $pIp$ have real rank zero. We do not need
separativity for $B$, just $K_1$-injectivity (by Lin's
result). Since
$I$ is an ideal of $A$, any projections in $M_\infty(I)$
which are (Murray-von Neumann) equivalent in $M_\infty(A)$
are also equivalent in
$M_\infty(I)$ (any implementing partial isometry necessarily
lies in
$M_\infty(I)$). Hence, the separativity of $A$ implies that
$I$ is separative, and so $pIp$ is separative. Therefore, by
Theorem 3.1, $pIp$ has $K_1$-surjectivity. 

With the above information in hand, Brown's proof \cite{\Bro,
Theorem 1} carries through in the present setting. We sketch
the details for the reader's convenience. By hypothesis,
$q=upu^*$ for some unitary $u\in U(B)$; let $\alpha$ denote
the image of $u$ in $K_1(B)$. Now $K_1(B)= K_1(I^\sim)=
K_1(I)$, and because $pIp$ is a full hereditary
sub-C*-algebra of $I$, the natural map $K_1(pIp) \rightarrow
K_1(I)$ is an isomorphism \cite{\Bro, Remark, p\. 117}. Thus
$\alpha$ is the image of some $\beta\in K_1(pIp)$. Since
$pIp$ has $K_1$-surjectivity,
$\beta$ is the image of some unitary $v_1 \in U(pIp)$. Let
$v= v_1 +1-p$ and $w=uv^*$. Then $w$ is a unitary in $B$ such
that $q=wpw^*$, and the image of $w$ in $K_1(B)$ is zero.
Since $B$ has
$K_1$-injectivity, $w\in U(B)^\circ$, from which it follows
that $p$ and
$q$ are homotopic. \qed\enddemo

\head Acknowledgements\endhead

We thank Michael Barr for a stimulating conversation and
Larry Levy for helpful comments.

\enddocument